\documentclass[14pt,leqno]{article}

\usepackage{graphicx}

\usepackage{palatino}
\usepackage[mathcal]{euler}

\usepackage{amsmath,amsthm,amsfonts,amssymb,latexsym,amscd,enumerate}

 \usepackage{comment}
\usepackage[small,nohug,heads=vee]{diagrams}
\diagramstyle[labelstyle=\scriptstyle]
\usepackage{xy}
\xyoption{all}
\theoremstyle{plain}
    \newtheorem{theorem}{Theorem}[section]
    \newtheorem{lemma}[theorem]{Lemma}
    
    \newtheorem{proposition}[theorem]{Proposition}

 \theoremstyle{definition}
    \newtheorem{definition}[theorem]{Definition}
    
    \newtheorem{example}[theorem]{Example}

    \newtheorem{remark}[theorem]{Remark}

\theoremstyle{remark}

\numberwithin{equation}{section}

\DeclareMathOperator{\Spin}{Spin}

\begin{document}


\newcommand{\myemph}{\emph}

\newcommand{\Spinc}{\Spin^c}

    \newcommand{\R}{\mathbb{R}}
    \newcommand{\C}{\mathbb{C}} 
    \newcommand{\N}{\mathbb{N}}
    \newcommand{\Z}{\mathbb{Z}} 
    \newcommand{\Q}{\mathbb{Q}}

\newcommand{\g}{\mathfrak{g}}
\newcommand{\h}{\mathfrak{h}}
\newcommand{\p}{\mathfrak{p}}
\newcommand{\kt}{\mathfrak{t}}
\newcommand{\M}{\mathcal{M}}

\newcommand{\cE}{\mathcal{E}}
\newcommand{\cS}{\mathcal{S}}
\newcommand{\cL}{\mathcal{L}}
\newcommand{\cH}{\mathcal{H}}
\newcommand{\cO}{\mathcal{O}}
\newcommand{\cP}{\mathcal{P}}
\newcommand{\cD}{\mathcal{D}}
\newcommand{\cF}{\mathcal{F}}
\newcommand{\cX}{\mathcal{X}}
\newcommand{\kA}{\mathfrak{A}}

\newcommand{\Sj}{ \sum_{j = 1}^{\dim M}}
\newcommand{\ii}{\sqrt{-1}}

\newcommand{\Wedge}{\Lambda}

\newcommand{\specialin}{\hspace{-1mm} \in \hspace{1mm} }

\title{Dirac operators on quasi-Hamiltonian $G$-spaces}  
\date{\today}
\author{Yanli Song\footnote{University of Toronto, \texttt{songyanl@math.utoronto.ca}} } 

\maketitle
\begin{abstract} 
We construct twisted spinor bundles as well as twisted pre-quantum bundles on quasi-Hamiltonian $G$-spaces, using the spin representation of loop group and the Hilbert space of Wess-Zumino-Witten model. We then define a Hilbert space together with a Dirac operator acting on it. The main result of this paper is that we show the Dirac operator has a well-defined index given by positive energy representation of loop group. This generalizes the geometric quantization of Hamiltonian $G$-spaces to quasi-Hamiltonian $G$-spaces. 
\end{abstract}
\tableofcontents

\abovedisplayskip=2pt
\belowdisplayskip=2pt

\section{Introduction}
Let $G$ be a compact, connected Lie group, and  $(M, \omega)$ a compact symplectic manifold with a Hamiltonian $G$-action. By choosing a $G$-invariant $\omega$-compatible almost complex structure on $M$, we can define a $G$-equivariant $\Z_2$-graded spinor bundle $S^\pm_{M}$. If the Hamiltonian $G$-space $M$ is pre-quantizable and has a $G$-equivariant pre-quantum line bundle $L$, we define a $\Z_{2}$-graded Hilbert space by 
\[
\mathcal{H}^\pm = L^{2}(M, S^\pm_{M}\otimes L)
\]
and a $G$-equivariant Spin$^{c}$-Dirac operator
\[
D^\pm : \cH^\pm \to \cH^\mp.
\]
Attributed to Bott, the \emph{quantization} of $(M, \omega)$ can be defined as the equivariant index 
\[
Q(M, \omega) = \mathrm{Ind}(D) =  [\mathrm{ker}(D^+)] - [\mathrm{ker}(D^-)] \in R(G).
\] 

The goal of this paper is to generalize the quantization process to the \emph{quasi-Hamiltonian $G$-space} introduced by Alekseev-Malkin-Meinrenken \cite{AMM98}. The q-Hamiltonian $G$-space, arising from infinite-dimensional  Hamiltonian loop group space, differs in many respects from Hamiltonian $G$-space. In particular, the moment map takes values in the group $G$ and the 2-form $\omega$ doesn't have to be closed or non-degenerate. Consequently,  the two key ingredients in defining $Q(M, \omega)$: the spinor bundle $S_{M}$ and pre-quantum line bundle $L$ might not exist in general. 

Given a q-Hamiltonian $G$-space $(M, \omega)$, we use the spin representation of loop group to construct twisted spinor bundles $S^\mathrm{spin}$ on $M$, and  the Hilbert space of Wess-Zumino-Witten model to construct twisted pre-quantum bundles $S^\mathrm{pre}$. Both of them are bundles of Hilbert spaces and play the same roles as the spinor bundle and pre-quantum line bundle for Hamiltonian $G$-spaces. We analogously define a Hilbert space
\[
\mathcal{H}:= \big[L^{2}(M, S^{\mathrm{spin}}\otimes S^{\mathrm{pre}})\big]^{G}. 
\]

One key in the construction of Dirac operators on $\cH$ is the algebraically defined cubic Dirac operator. It was introduced by Kostant for finite-dimensional Lie group, and extensively studied for infinite-dimensional loop group by various people.  Our strategy is to construct a Dirac operator as a combination of algebraic cubic Dirac operators and geometric Spin$^{c}$-Dirac operators. To be more precise, we choose an open cover of $M$ using the symplectic cross-section theorem for q-Hamiltonian $G$-spaces,  so that every open subset $U$ has the geometric structure:
\[
U \cong G \times_{H} V,
\]
where the slice $V$ is a Hamiltonian $H$-space. Accordingly, the tangent bundle $TU$ splits equivariantly into ``vertical direction" and ``horizontal direction". We define a suitable Dirac operator on $U$ so that it acts as the Spin$^{c}$-Dirac operator on the vertical part $V$ and the cubic Dirac operator for loop group on the horizontal part $G/H$. Using partition of unity, we obtain a global Dirac operator $\cD$ on $\cH$ by patching together Dirac operators on the open sets $U$.  The main result of this paper is that we show the Dirac operator $\cD$ has a well-defined index given by positive energy representations of loop group.

\subsection*{Acknowledgements}
The author received many useful suggestions as well as warm help from Eckhard Meinrenken. The author also would like to thank Nigel Higson and Mathai Varghese  for inspiring conversations.

\section{Loop group and positive energy representation}
We first give a brief review on loop groups and their representations. We use \cite{PS86} as our primary reference. 

\subsection{Loop group and central extension}
Let $G$ be a compact, simple and simply connected Lie group, and fix a ``Sobolev level" $s > 1$. We define $LG$ the \emph{loop group}  as the Banach Lie group consisting of maps $S^{1} \to G$ of Sobolev class $s +\frac{1}{2}$ with the group structure given by pointwise multiplication. The Lie algebra $L\g = \Omega^{0}(S^{1}, \g)$ is given by the space Lie algebra $\g$-valued 0-forms of Sobolev class $s + \frac{1}{2}$ and  $L\g^{*} = \Omega^{1}(S^{1}, \g)$  the space of $\g$-valued 1-forms of Sobolev class $s-\frac{1}{2}$. Integration over $S^{1}$ gives a natural non-degenerate pairing between $L\g$ with $L\g^{*}$.

Note that $L\g^{*}$ can be identified with the affine space of connections on the trivial principle $G$-bundle over $S^{1}$. The loop group $LG$ acts on $L\g^{*}$ by gauge transformation
\begin{equation}
\label{gauge transformation}
g \cdot \xi = \mathrm{Ad}_{g}(\xi) - dg \cdot g^{-1}, \hspace{5mm} g \in LG, \xi \in L\g^{*},
\end{equation}
where $dg\cdot g^{-1}$ is the pull-back of the right-invariant Maurer-Cartan form on $G$. 

Let $\widehat{LG}$ be the basic central extension of $LG$, defined infinitesimally by the cocycle
\[
(\xi_1, \xi_2) \mapsto \oint d\xi_1 \cdot \xi_2, \hspace{5mm} \xi_1,\xi_2 \in L\g^*. 
\]
The coadjoint action of $LG$ on 
\[
\widehat{L\g}^* = L\g^*\oplus \R
\] 
is given by the formula
\[
g \cdot (\xi,k) = (\mathrm{Ad}_{g}(\xi) - k \cdot g^{-1}dg, k). 
\]
One can view the action (\ref{gauge transformation})  as the coadjoint action on the affine hyperplane $L\g^* \times \{1\} \subset \widehat{L\g}^*$.

Fixing a maximal torus $T$, the choice of a set of positive roots $\mathfrak{R}_{+}$ for $G$ determines a positive Weyl chamber $\mathfrak{t}_{+}^*$. It is well-known that the orbits of coadjoint $G$-action on $\g^{*}$ are parametrized by points in $\kt^*_+$. The set of coadjoint $LG$-orbits can be described as follow. Denote by $\alpha_{0}$ the highest root and 
\[
\rho_{G} = \frac{1}{2}\sum_{\alpha \in \mathfrak{R}_{+}}\alpha.
\]
There is a unique ad-invariant inner product $\langle \cdot, \cdot \rangle_{\g}$ on $\g$, rescaled so that the highest root  of $\g$ has norm $\sqrt{2}$. The \emph{dual Coxeter number} of $G$ is defined by 
\[
h^{\vee} = 1 + \langle \rho_{G}, \alpha_{0} \rangle_\g,
\]
and the \emph{fundamental Weyl alcove} for $G$ is the simplex
\[
\mathfrak{A} = \{ \xi \in \mathfrak{t}_{+}\big| \langle \alpha_{0}, \xi \rangle_{\g} \leq 1\} \subset \mathfrak{t} \subset \g. 
\]
Every coadjoint orbit of $LG$-action on $L\g^*$ contains a unique point in $\kA$.

For any $\xi \in L\g^*$, we define the \emph{holonomy map} 
\[
\mathrm{Hol}: L\g^* \to G
\]
the smooth map that sends $\xi$ to the holonomy of $\xi$ around $S^1$. This map satisfies the equivariance property 
\[
\mathrm{Hol}(g \cdot \xi) = g(0) \cdot \mathrm{Hol}(\xi) \cdot g(0)^{-1}. 
\]
It follows that the based loop group
\[
\Omega G = \{g \in LG \big| g(0) = e\}
\]
acts freely on $L\g^*$. 

The isotropy group $(LG)_{\xi}$ is isomorphic to $G_{\mathrm{Hol}(\xi)}$, and thus compact. It sets up a 1-1 correspondence between the set of coadjoint $LG$-orbits and conjugacy classes. Moreover $(LG)_{\xi}$ with $\xi \in \kt^*_+$ depends only on the open face $\sigma$ of $\kA$ containing $\xi$ and will be denoted by $(LG)_{\sigma}$. If we  introduce a partial order on open faces by setting $\tau \preceq\sigma$ if $\tau \subseteq \overline{\sigma}$, then one has that
\[
\sigma \preceq \tau \Rightarrow (LG)_{\tau} \subseteq (LG)_{\sigma}.\]
In particular, $(LG)_{0} = G$ and $(LG)_{\mathrm{int} \mathfrak{A}} = T$.

\subsection{Positive Energy Representation}
Let $S^{1}_{\mathrm{rot}}$ be the rotation group on $S^{1}$ and $\partial$ its infinitesimal generator. 
Consider a unitary representation of $S^{1}_{\mathrm{rot}}\ltimes \widehat{LG}$ on a Hilbert space $V$, on which the central circle acts by scalar multiplication. 

\begin{definition}
We say that $V$ is a \emph{positive energy representation} if  the  operator $\partial$ is self-adjoint with spectrum bounded below. Moreover, we say that $V$ \emph{has level k} if the central circle of $\widehat{LG}$ acts with weight $k$. 
\end{definition} 

The positive energy representations of loop groups behave quite analogously to the representation theory of compact Lie groups. For example, every irreducible positive energy representation is uniquely determined by the highest weight. To be more precise, let $T$  be a maximal torus of $G$ and $\Lambda^{*}$ the weight lattice. We take $S^{1}_{\mathrm{rot}} \times T \times S^{1}$ as the maximal torus of $S^{1}_{\mathrm{rot}}\ltimes \widehat{LG} $, where the second $S^{1}$ factor comes from the central extension. The \emph{affine weights} of $LG$ are in the forms of $(m, \lambda, k)$, where $m \in \Z$ is the \emph{energy}, $\lambda \in \Lambda^{*}$ is the weight of $G$, and $k$ is the \emph{level}. 

The \emph{affine Weyl group}  $W_{\mathrm{aff}} = W \ltimes \Lambda$ acts on affine weights as follow: the Weyl group $W$ acts as usual on $\Lambda^{*}$ and the action of $z \in \Lambda$ is given by 
\[
z \cdot (m,\lambda, k) = (m + \langle \lambda, z \rangle + \frac{k}{2} \cdot \|z\|^{2}, \lambda+ k \cdot z, k). 
\]
The level $k$ is fixed by the affine Weyl group action and the energy is shifted so as to preserve the inner product:
\begin{equation}
\label{weight product}
(m_{1}, \lambda_{1}, k_{1}) \cdot (m_{2}, \lambda_{2}, k_{2})= \langle \lambda_{1}, \lambda_{2} \rangle - m_{1}k_{2} - m_{2}k_{1}. 
\end{equation}
For a fixed level $k$, every irreducible positive energy representation $V$ of $LG$ is uniquely determined by the dominant weight $\lambda$ at the minimum energy $m$. We call $\boldsymbol{\lambda}  = (m, \lambda, k)$ the \emph{highest weight} of $V$. 

Let $V(n)$ be a subspace of the Hilbert space $V$, on which the energy operator acts on $V(n)$ with weight $n$. The positive energy condition asserts that there is an integer $n_{\mathrm{min}}$ so that $V(n) = 0$ for all $n < n_{\mathrm{min}}$. The algebraic direct sum
\[
V^{\mathrm{fin}}:= \bigoplus_n V(n)
\] 
consists of vectors of \emph{finite energy}. It is a dense subspace of $V$. In addition, we can always normalize so that the lowest energy level $n_{\mathrm{min}}$ equals zero.

It is well-known by the Borel-Weil theorem that  all the irreducible $G$-representations are parameterized by the dominant weights $P_{G,+} = \Lambda^{*} \cap \mathfrak{t}^{*}_{+}$. Similarly there is a 1-1 correspondence between irreducible positive energy representation at level $k$ and weight in 
\[
P_{k,+} = k \mathfrak{A} \cap \Lambda^{*} =  \{ \lambda \in \Lambda^{*} \big| \frac{\lambda}{k} \in \mathfrak{A} \}.
\]
The abelian group $R_{k}(LG)$ generated by irreducible positive energy representations at level $k$  has a finite basis and a ring structure known as \emph{fusion product}. But we won't discuss it in this paper.

\section{Dirac operators in the algebraic setting}
The \emph{cubic Dirac operator} is an algebraically defined operator introduced by Kostant \cite{kostant99} for finite-dimensional Lie algebras. It has now been generalized to infinite-dimensional case and plays an important role in the theory of loop groups: its application in representation theory was first demonstrated in the lecture notes of Wassermann \cite{Wassermann-notes}. Later Landweber and Posthuma generalize it to different homogeneous settings \cite{landweber01, Posthuma11}. A family version of the cubic Dirac operator was  used by Freed-Hopkins-Teleman \cite{Freed13} to construct the isomorphism between the twisted K-theory and fusion ring of loop groups. In addition Meinrenken \cite{MR2894443} discusses its application in general Kac-Moody algebra. 

\subsection{Finite dimensional case}
Let $G$ be a compact Lie group and $\g$ its Lie algebra equipped with an ad-invariant inner product $\langle \ , \ \rangle_\g$. Let $\mathrm{Cliff}(\g)$ be the $\Z_2$-graded complex Clifford algebra of $\g$ and  $S_{\g}$ an irreducible $\Z_{2}$-graded $\mathrm{Cliff}(\g)$-module. 

Fix an orthonormal basis 
\[
X_{a}, \hspace{5mm} a = 1, \dots, \mathrm{dim}\g.
\]
We define a map $\mathrm{ad}^{\g}: \g \to \mathrm{Cliff}(\g)$ by the formula 
\[
\mathrm{ad}^{\g}(X) := \frac{1}{4}\cdot \sum_{a=1}^{\mathrm{dim}\g} [X, X_{a}] \cdot X_{a}, \hspace{5mm} X \in \g.
\]

Let $U(\g)$ be the universal enveloping algebra of $\g$. We consider the \emph{non-commutative Weil algebra} 
\[
U(\g) \otimes \mathrm{Cliff}(\g) 
\]
introduced by Alekseev-Meinrenken\cite{Meinrenken00}.
\begin{definition}
\label{cubic operator}
The \emph{cubic Dirac operator} $D_\g$ is an element of the 
algebra $U(\g) \otimes \mathrm{Cliff}(\g)$ defined as
\[
D_{\g}= \sum_{a=1}^{\mathrm{dim}\g} \big( X_{a} \otimes X_{a} +  \frac{1}{3} \otimes  \mathrm{ad}^{\g}(X_{a}) \cdot X_{a} \big)
\]
\end{definition}
The key property of the cubic Dirac operator is that its square is simple and nice.

More generally, let $H \subset G$ be a closed subgroup of the equal rank. Using the inner product we write 
\[
\mathfrak{g} = \mathfrak{h} \oplus \mathfrak{p}, \hspace{5mm} \p = \h^\perp.
\] 
This decomposition induces isomorphisms:
\[
\mathrm{Cliff}({\g}) \cong \mathrm{Cliff}({\h}) \otimes \mathrm{Cliff}({\mathfrak{p}}), \hspace{5mm}
S_{\g} \cong S_{\h} \otimes S_{\mathfrak{p}},
\]
where $S_{\h}, S_{\p}$ are spinor modules of $\mathrm{Cliff}({\h})$ and $\mathrm{Cliff}({\p})$ respectively. 

\begin{definition}
We define the \emph{relative cubic Dirac operator} 
\[
D_{\g,\h} \in U(\g) \otimes \mathrm{Cliff}(\p)
\] 
by 
\begin{equation}
\label{relative cubic operator}
D_{\g, \h} =\sum^{(\p)}\big(X_{a} \otimes X_{a} + \frac{1}{3} \otimes   \mathrm{ad}^{\p}(X_{a}) \cdot X_{a} \big).
\end{equation}
Here $\sum^{(\p)}$ indicates the summation over basis of $\p$. 

\end{definition}

As it stands, the element $D_{\g, \h}$ gives us an operator on  $W \otimes S_{\p}$ for any $\g$-representation $W$. To exhibit the structure of $D_{\g,\h}$, we decompose $W\otimes S_{\p}$ with respect to the $\h$-action and denote by $M(\nu)$ the isotypic $\h$-summand with highest weight $\nu$.

\begin{theorem}[\cite{kostant99}]
\label{relative cubic}
Suppose that $W_{\lambda}$ is an irreducible $\g$-representation with highest weight $\lambda$. The following formula holds
\[
D_{\g,\h}^{2}\big|_{M(\nu)} = \|\lambda +\rho_{G} \|^{2}- \|\nu + \rho_{H}\|^{2}. 
\]
\end{theorem}

\subsection{Dirac operators on Homogeneous spaces}

Suppose now that $M = G/H$ is an orbit of the coadjoint $G$-action on $\g^*$. It is known that $M$ has a $G$-invariant complex structure, which determines a $H$-equivariant splitting: $\p = \p^+ \oplus \p^-$. One can check that the spinor bundle associated to the complex structure on $M$ is given by 
\begin{equation}\label{spinor for lg}
S_{M} = G \times_H (\wedge \p^-) \cong G \times_{H} (S_{\p}^{*} \otimes \C_{\rho_G - \rho_H}),
\end{equation}
and the canonical line bundle
\begin{equation}\label{canonical for lg}
K_M = G \times_H\C_{2 \cdot(\rho_G - \rho_H)}. 
\end{equation}
Hence, the Hilbert space $L^{2}(M, S_{M})$ can be identified with 
\begin{equation}\label{G dec}
 \big[ L^{2}(G) \otimes S_{\p}^{*} \otimes  \C_{\rho_G - \rho_H}\big]^{H}  \cong \bigoplus_{\lambda \in P_{G,+}} W_{\lambda} \otimes [W_{\lambda}^{*} \otimes S_{\p}^{*} \otimes  \C_{\rho_G - \rho_H}]^{H},
\end{equation}
where the isomorphism comes from the Peter-Weyl theorem. 

We  define Dirac operators on $M$ in two different ways. First of all, the Levi-Civita connection $\nabla^{TM}$ lifts to a Hermitian connection $\nabla^{S_{M}}$ on $S_{M}$. In particular, the connection $\nabla^{S_{M}}$ is compatible with the Clifford action in the sense that
\[
[X, \nabla^{S_{M}}_Y] = \nabla^{TM}_Y X, \hspace{5mm} X, Y \in TM.  
\]
We define a geometric Spin$^{c}$-Dirac operator by
\[
D_{\mathrm{geo}} = \sum_{a=1}^{\mathrm{dim}M} X_{a} \cdot \nabla^{S_{M}}_{X_{a}},
\]
where $\{X_{a}\}_{i=1}^{\mathrm{dim}M}$ is an orthonormal basis of $TM$.  

On the other hand, since the cubic Dirac operator $D_{\g,\h}$ is $H$-equivariant, it restricts to an operator on 
\[
[W_{\lambda}^{*} \otimes S_{\p}^{*}  \otimes  \C_{\rho_G - \rho_H}]^{H}.
\] 
Tensoring  the identity operator on each $W_{\lambda}$, and summing over $W_{\lambda}$, one obtains an  operator $D_{\mathrm{alg}}$ on (\ref{G dec}). 

\begin{lemma}
\label{geo alg}
The difference between $D_{\mathrm{alg}}$ and $D_{\mathrm{geo}}$ on (\ref{G dec}) is a bounded operator. 
\begin{proof}
We rewrite the geometric connection $\nabla^{S_M}$  as 
\[
\nabla^{S_M}_X = X + \frac{1}{2}\mathrm{ad}^{\p}(X). 
\]
It follows that
\[
D_{\mathrm{geo}} =\sum_{i=1}^{\mathrm{dim}\p} X_{a} \otimes \big(X_{a} + \frac{1}{2}\mathrm{ad}^{\p}(X_{a})\big).
\]
We deduce the lemma by comparing the above with (\ref{relative cubic operator}).  An alternative proof can be found in \cite[Chapter 9]{Meinrenken-book}.
\end{proof}
\end{lemma}

\subsection{Infinite dimensional case}
The definitions of spin representation and cubic Dirac operator can be extended to the infinite-dimensional loop algebra $L\g$. 

Let now $G$ be a compact, simple and simple connected Lie group with Lie algebra $\g$. The loop algebra $L\g$ carries an inner product defined by  
\begin{equation}
\label{inner product}
B(X, Y) = \frac{1}{2\pi}\int_{0}^{2\pi} \langle X(\theta), Y(\theta) \rangle_{\g} d\theta, \hspace{5mm} X, Y \in L\g. 
\end{equation}
As in the finite dimensional case, we can define the Clifford algebra $\mathrm{Cliff}(L\g)$, and its spin representation $S_{L\g}$. Here $S_{L\g}$ is a $\Z_2$-graded complex Hilbert space, and also a positive energy $LG$-representation with highest weight 
\[
\boldsymbol{\rho}_G= (0, \rho_G, h^{\vee}).
\]
The explicit construction of $S_{L\g}$ is given in \cite{PS86}. For the general theory of Clifford algebras and  representations for infinite dimensional Hilbert spaces we refer to \cite{Plymen94}.

Let us fix an orthonormal basis $\{X_{a}\}$  of $\g$. For $n \in \Z$, we write $X_{a}^{n}$ for the loop 
\[
s \mapsto e^{i n s}\cdot X_{a}, \hspace{5mm} s \in \R,
\] 
and $\g(n)$ the vector space spanned by $\{X^{n}_{a}\}_{a=1}^{\mathrm{dim}\g}$. 
The algebraic direct sum
\[
L\g^{\mathrm{fin}}:= \bigoplus_{n \in \Z} \g(n)
\] 
is dense in $L\g_{\C}$. The dense subspace of the spin representation $S_{L\g}$ which consists of vectors with  finite energy may be realized as
\[
S_{L\g}^{\mathrm{fin}}= S_{\g} \otimes \bigotimes_{k>0} \Lambda^*\big(\g_\C z^k\big).
\]

\begin{definition}
We define the \emph{cubic Dirac operator} 
\[
D_{L\g} \in U(L\g) \otimes \mathrm{Cliff}(L\g)
\] 
by
\begin{equation}\label{d_lg}
D_{L\g} = \sum_{n\in \Z} \sum_{a=1}^{\mathrm{dim}\g} \big( X_{a}^{n}\otimes X_{a}^{-n}+ \frac{1}{3}\otimes \mathrm{ad}^{L\g}(X_{a}^{n}) \cdot X_{a}^{-n} \big),
\end{equation}
\end{definition}
\begin{remark}
The map
\begin{equation}\label{ad_lg}
\mathrm{ad}^{L\g}(X) := \frac{1}{4} \sum_{n\in \Z} \sum_{a=1}^{\mathrm{dim}\g} [X, X_{a}^n] \cdot X_{a}^{-n} \in \mathrm{Cliff}(L\g), \hspace{5mm} X \in L\g^{\mathrm{fin}}.
\end{equation}
is defined only on vectors with finite energy. To justify the infinite summation in (\ref{d_lg}) and (\ref{ad_lg}), we refer to \cite{MR2894443, Wassermann-notes}. For any positive energy $LG$-representation $V$, $D_{L\g}$ gives us an unbounded operator on $\big(V \otimes S_{L\g}\big)^{\mathrm{fin}}$ which is a  dense subspace of $V \otimes S_{L\g}$.
\end{remark}

Let $\h$ be an isotropy Lie algebra $\h$ of the coadjoint $LG$-action. We decompose the Lie algebra $L\g = \p \oplus \h$, where $\p = \h^{\perp}$. By the multiplicative property of spin representation 
\[
S_{L\g} = S_{\p} \otimes S_{\h},
\]
where $S_{\p}$  is a $\Z_2$-graded irreducible representation of the Clifford algebra $\mathrm{Cliff}(\p)$. 

\begin{definition}
We define  $D_{L\g, \h}\in U(L\g) \otimes \mathrm{Cliff}(\p)$ by the formula
\[
D_{L\g, \h} = \sum_{a, n}^{(\p)}\big( X_a^n \otimes X_a^{-n} + \frac{1}{3} \otimes \mathrm{ad}^\p(X_{a}^n)\cdot X_a^{-n} \big),
\]
where the summation ranges over a basis of $\p$. 
\end{definition}

We denote by $\widehat{H}$ the central extension of $H$ induced by the inclusion $H \hookrightarrow LG$ and the central extension $\widehat{LG}$. We decompose $V \otimes S_{\p}$ with respect to the $S^{1}_{\mathrm{rot}} \times \widehat{H}$-action and denote by $M(\boldsymbol{\nu})$ the isotypic  component labeled by 
\[
\boldsymbol{\nu} = (n, \nu, k+h^\vee) \in \Z \times \Lambda^{*} \times \Z.
\] 
We have an analog of Theorem \ref{relative cubic} for the infinite-dimensional case. 
\begin{theorem}
\label{square operator-4}
Suppose that $V_{\lambda}$ is an irreducible positive energy representation with highest weight $\boldsymbol{\lambda}  = (0, \lambda, k)$. If we restrict to the isotypic component $M(\boldsymbol{\nu})$ of $V_{\lambda} \otimes S_{\p}$, we have that
\[
D_{L\g,\h}^{2}\big|_{M(\boldsymbol{\nu})}= \|\boldsymbol{\lambda} + \boldsymbol{\rho}_{G}\|^{2} - \|\boldsymbol{\nu} + \boldsymbol{\rho}_{H}\|^{2},
\]
where $\boldsymbol{\rho}_{G} = (0, \rho_{G}, h^{\vee})$ and $ \boldsymbol{\rho}_{H} = (0, \rho_{H}, 0)$.
\begin{proof}
\cite[Theorem 7.5]{MR2894443}.
\end{proof}
\end{theorem}

\section{Hamiltonian $LG$-spaces and q-Hamiltonian $G$-spaces}
The theory of q-Hamiltonian $G$-spaces was developed in \cite{AMM98}. It provides a finite-dimensional model for Hamiltonian $LG$-spaces. In this section, we begin by reviewing the basic definitions, and then discuss their cross-section theorems. We assume that $G$ is a compact, simple and simple connected Lie group. 

\subsection{Basic definitions}
Recall that a \emph{Hamiltonian $G$-space} is a triple $(M, \omega, \mu)$, with $\omega$ the $G$-equivariant symplectic 2-form, and $\mu: M \to \g^{*}$ the \emph{moment map} satisfying that 
\[
\iota_{\xi_{M}} \omega =d \langle \mu, \xi \rangle, \hspace{5mm} \xi \in \g, 
\]
where $\xi^M$ is the vector field on $M$ induced by the infinitesimal action of $\xi$. 

The above definition can be extended to the loop group setting. Let $\mathcal{M}$ be an infinite-dimensional Banach manifold. We say that it is \emph{weakly symplectic} if it is equipped with a closed 2-form $\omega \in \Omega^{2}(\mathcal{M})$ so that the induced map
\[
\omega^{\flat}: T_{m}\mathcal{M} \to T_{m}^{*}\mathcal{M}
\]
is injective. 

\begin{definition}
A \emph{Hamiltonian $LG$-space} is a weakly symplectic Banach manifold $(\mathcal{M}, \omega)$ together with a $LG$-action and 
a $LG$-equivariant map $\mu: \mathcal{M} \to L\g^{*}$ so that $\iota_{\xi_{\mathcal{M}}} \omega =d \langle \mu, \xi \rangle$ for all $\xi \in L\g$.
\end{definition}
For example, the coadjoint $LG$-orbit is a Hamiltonian $LG$-space, with moment map the inclusion. 

Let $(\mathcal{M}, \omega, \mu)$ be a Hamiltonian $LG$-space. Since the based loop group $\Omega G$ acts freely on $L\g^{*}$, it acts freely on $\mathcal{M}$ as well by the equivariance of $\mu$. We thus obtain a commuting square
\begin{diagram}
 \mathcal{M} & \rTo^{\mu} &   L\g^{*}\\
\dTo &  &\dTo^{\mathrm{Hol}}\\
M & \rTo^{\phi} & G\\,
\end{diagram}
where the quotient $M = \mathcal{M}/\Omega G$ is a finite-dimensional compact smooth manifold provided that $\mu$ is proper. 

Alekseev-Malkin-Meinrenken \cite{AMM98} give a set of conditions a $G$-space $M$ must satisfy in order to arise from a Hamiltonian $LG$-space by such a construction.  

Choose an invariant inner product $\langle \cdot, \cdot \rangle_{\g}$ on $\g$ and denote by $\theta^{L}, \theta^{R} \in \Omega^{1}(G,\g)$ the left and right invariant Maurer-Cartan forms on $G$ and  the \emph{Cartan 3-form}  
\[
\chi = \frac{1}{12} \langle \theta^{L}, [\theta^{L}, \theta^{L}] \rangle_{\g} = \frac{1}{12} \langle \theta^{R}, [\theta^{R}, \theta^{R}]\rangle_{\g} \in \Omega^{3}(G). 
\]

\begin{definition}[\cite{AMM98}]
\label{q-Hamiltonian}
A \emph{q-Hamiltonian} $G$-space is a compact $G$-manifold $M$, together with an equivariant 2-form $\omega$, and an equivariant  map $\phi: M \to G$ satisfying the following properties:
\begin{enumerate}
\item
$d \omega = \phi^{*}\chi;$
\item
$\iota_{\xi_{M}} \omega = \frac{1}{2}\langle \phi^{*}(\theta^{L} + \theta^{R}), \xi \rangle_{\g}$ for all  $\xi \in \g;$
\item
$\mathrm{ker}(\omega) \cap \mathrm{ker}(d\phi) = 0. $
\end{enumerate}
We call $\phi$ the \emph{group-valued moment map}. 
\end{definition}
According to \cite[Theorem 8.3]{AMM98}, there is a 1-1 correspondence between Hamiltonian $LG$-spaces with proper moment map and q-Hamiltonian $G$-spaces. One can always choose to work with infinite-dimensional Hamiltonian $LG$-spaces with more conventional definitions or to use finite dimensional q-Hamiltonian $G$-spaces. The counterparts of coadjoint orbits for q-Hamiltonian $G$-spaces are conjugacy classes $\mathcal{C}$ in $G$ with group-valued moment map the embedding $\mathcal{C} \hookrightarrow G$.

\subsection{Cross-section theorems}
The Hamiltonian $LG$-spaces and their equivalent finite-dimensional models behave in many respects like the usual Hamiltonian $G$-spaces. This is due to the existence of the cross-section theorem we shall now describe.  

 Let us first introduce a partial order of open faces of $\mathfrak{A}$ by setting $\tau \preceq\sigma$ if $\tau \subseteq \overline{\sigma}$. The isotropy group $(LG)_{\xi}$ of the coadjoint $LG$-action on $L\g^{*}$ depends only on the open face $\sigma$ of $\mathfrak{A}$ containing $\xi$ and will be denoted by $(LG)_{\sigma}$ (note however $(LG)_{\sigma}$ will generally contain non-constant loops). One has that
\[
\sigma \preceq \tau \Rightarrow (LG)_{\tau} \subseteq (LG)_{\sigma}.
\]
In particular, $(LG)_{0} = G$ and $(LG)_{\mathrm{int} \mathfrak{A}} = T$. 

We define a $(LG)_{\sigma}$-invariant open subset of $(L\g)_{\sigma}^{*}$ by
\[
\mathcal{A}_{\sigma} = (LG)_{\sigma} \cdot \bigcup_{\sigma \preceq\tau} \tau. 
\]
Note that $\mathcal{A}_{\sigma}$ is a slice for all $\xi \in \sigma$ for the action of $LG$ in the sense that
\[
LG \times_{(LG)_{\sigma}} \mathcal{A}_{\sigma} \to LG \cdot \mathcal{A}_{\sigma}
\]
is a diffeomorphism of Banach manifolds. 
\begin{theorem}
\label{cross section LG}
Let $(\mathcal{M}, \omega, \mu)$ be a Hamiltonian $LG$-space with proper moment map. For every open face $\sigma$ of $\mathfrak{A}$, the cross-section
\[
\mathcal{V}_{\sigma} = \mu^{-1}(\mathcal{A}_{\sigma})
\]
is a finite-dimensional symplectic submanifold with Hamiltonian $(LG)_{\sigma}$-action. The restriction of $\mu|_{\mathcal{V}_{\sigma}}$ is a moment map of the $\widehat{(LG)}_{\sigma}$-action (the central circle acts trivially on $\mathcal{V}_{\sigma}$).
\begin{proof}
\cite[Theorem 4.8]{MW98}.
\end{proof}
\end{theorem}

The symplectic cross-section theorem  carries over to q-Hamiltonian $G$-spaces. The centralizer $G_{\mathrm{exp}(\xi)}$ with $\xi \in \mathfrak{A}$ is isomorphic to $(LG)_{\xi}$ and it depends only on the open face $\sigma$ of $\mathfrak{A}$ containing $\xi$. We denoted it by $G_{\sigma}$. The subset
\[
 A_{\sigma} = \mathrm{Ad}(G_{\sigma}) \cdot \mathrm{exp}(\bigcup_{\sigma \preceq\tau}\tau) \subset G_{\sigma} \subset G
\]
is smooth and is a slice for the $\mathrm{Ad}$(G)-action at points in $\sigma$. 
\begin{theorem}
 Let $(M, \omega, \phi)$ be a q-Hamiltonian $G$-space. The cross-section
\begin{equation}
\label{slice}
V_{\sigma} = \phi^{-1}(A_{\sigma}) 
\end{equation}
 is a smooth $G_{\sigma}$-invariant submanifold and 
\[
G \times_{G_{\sigma}} V_{\sigma} \cong G \cdot V_{\sigma}
\]
is a $G$-invariant open subset of $M$. Moreover, $V_{\sigma}$
 is a q-Hamiltonian $G_{\sigma}$-space with the restriction of $\phi$ as the group-valued moment map.
\begin{proof}
\cite[Proposition 7.1]{AMM98}.
\end{proof}
\end{theorem}

\begin{remark}
\label{identify}
It is important to point out that if we identify $G_{\sigma}\cong (LG)_{\sigma}$, the two cross-sections
\[
 \mathcal{V}_{\sigma}\subset \mathcal{M}, \hspace{5mm} V_{\sigma} \subset M
\] 
are equivariantly diffeomorphic. In particular, every $V_\sigma$ is also a Hamiltonian $G_\sigma$-space. 
\end{remark}

\section{Twisted spinor bundle and twisted pre-quantum bundle}
In this section, we construct twisted spinor bundle and twisted pre-quantum bundles on q-Hamiltonian $G$-spaces.

\subsection{Construction of the twisted spinor bundle}
Let $G$ be a compact, simple, and simply connected Lie group. Let $(M, \omega)$ be a q-Hamiltonian $G$-space and $\mathcal{M}$ its corresponding Hamiltonian $LG$-space. 

We first replace the cross-sections $V_\sigma$ in (\ref{slice}) with smaller open subsets. To be more precise, for every vertex $\sigma$ of $\mathfrak{A}$, let  $Y_{\sigma}$ be a  $G_{\sigma}$-invariant,  open subset of $V_\sigma$ so that $\overline{Y_\sigma} \subset V_\sigma$, and
\[
M/G \subseteq \bigcup_{\sigma, \mathrm{dim}\sigma =0}  Y_{\sigma}.
\]
Then we form an open cover of $M$ by 
\[
\{ U_{\sigma} = G \times_{G_{\sigma}}Y_{\sigma}\}_{\sigma, \mathrm{dim}\sigma = 0}.
\] 
For all open faces $\tau$ of $\mathfrak{A}$ with $\mathrm{dim}\tau > 0$, we define
\[
Y_{\tau} = \bigcap_{\sigma \preceq\tau, \mathrm{dim}\sigma = 0} Y_{\sigma},  \hspace{5mm} U_\tau = G \times_{G_\tau} Y_{\tau}. 
\]
Remark at this point that each $Y_{\tau}$ is a Hamiltonian 
$G_\tau$-space and admits $G_\tau$-invariant almost complex structures.

\begin{lemma}
\label{homotopy}
There exists a collection of $G_{\sigma}$-invariant almost complex structures on the collection of $Y_{\sigma}$ such that the embedding 
\[
Y_{\tau} \hookrightarrow Y_{\sigma},\hspace{5mm} \sigma \prec \tau
\] 
is almost complex. In addition, any two almost complex structures with the required properties are homotopic. We denote by $S_{Y_{\sigma}}$ the spinor bundle on $Y_{\sigma}$ associated to the almost complex structures. 
\begin{proof}
\cite[Lemma 3.2]{Meinrenken-W-2001}.
\end{proof}
\end{lemma}

Let $\pi : U_\sigma \to G/G_\sigma$ be the projection. The tangent bundle decomposes $G$-equivariantly
\[
TU_\sigma \cong \pi^*T(G/G_\sigma) \oplus G \times_{G_\sigma} TY_\sigma. 
\]
The base manifold $G/G_\sigma$ is a conjugacy class and might not have a $G$-equivariant Spin$^c$-structure in general. Thus the total space $U_\sigma$ doesn't have to be Spin$^c$ either. 
 
On the other hand, the coadjoint $LG$-orbit $\mathcal{O} =LG/(LG)_\sigma$ is a complex manifold. By the discussion in  \cite{Meinrenken-W-2001}, the weight 
\[
2 \cdot (\rho_{G} -\rho_{\sigma}, h^{\vee})
\] 
is fixed by $\widehat{(LG)}_\sigma$ and the tensor product 
\[
S^{*}_{L\g/(L\g)_{\sigma}} \otimes \C_{(\rho_{G} -\rho_{\sigma}, h^{\vee})}
\] 
is a $(LG)_\sigma$-space. The associated spinor bundle and canonical line bundle on $\cO$ are given  by 
\[
S_\mathcal{O} = \widehat{LG} \times_{\widehat{(LG)}_\sigma} (S^{*}_{L\g/(L\g)_{\sigma}} \otimes \C_{(\rho_{G} -\rho_{\sigma}, h^{\vee})})
\]
and 
\[
K_\cO = \widehat{LG} \times_{\widehat{(LG)}_\sigma}  \C_{2 \cdot (\rho_{G} -\rho_{\sigma}, h^{\vee})}).
\]
 One can compare them with (\ref{spinor for lg}) and (\ref{canonical for lg}).

Motivated by the above, we define a bundle of Hilbert space on $U_\sigma$ by 
\begin{equation}
\label{spinor-2}
S_{U_{\sigma}}^{\mathrm{spin}} = G \times_{G_{\sigma}} \big(S^{*}_{L\g/(L\g)_{\sigma}} \otimes \C_{(\rho_{G} -\rho_{\sigma}, h^{\vee})} \otimes S_{Y_{\sigma}}\big),
\end{equation}
where $G_{\sigma}$ acts on 
\[
S^{*}_{L\g/(L\g)_{\sigma}} \otimes \C_{(\rho_{G} -\rho_{\sigma}, h^{\vee})} 
\] 
factoring through the identification $G_{\sigma} \cong (LG)_{\sigma}$. In addition, we equip $S_{U_{\sigma}}^{\mathrm{spin}}$ with a $\Z_{2}$-grading induced by that on $S^{*}_{L\g/(L\g)_{\sigma}}$ and $S_{Y_\sigma}$. To sum up, we obtain a collection of $G$-equivariant bundles of $\Z_{2}$-graded Hilbert space $\{ S_{U_{\sigma}}^{\mathrm{spin}}\}$.  We next show that such a collection of bundles can be glued together.
\begin{lemma}
For any $\sigma \prec \tau$, the normal bundle $\nu_{\tau}^{\sigma}$ of $Y_{\tau} \hookrightarrow Y_{\sigma}$ has a $G_{\tau}$-equivariant almost complex structure with spinor bundle isomorphic to
\begin{equation}
\label{normal} S_{(L\g)_{\sigma}/(L\g)_{\tau}}^{*} \otimes \C_{(\rho_{\sigma} - \rho_{\tau})},
\end{equation}
where $\rho_{\tau}, \rho_{\sigma}$ are the half-sums of positive roots for $G_\tau, G_\sigma$ respectively. 
\begin{proof}
By the cross-section theorem, the normal bundle $\nu_{\tau}^{\sigma}$ is  isomorphic to the trivial bundle 
\[
\g_{\sigma}/\g_{\tau} \cong (L\g)_{\sigma}/(L\g)_{\tau}
\] 
with equivariant almost complex structure. In fact,  let  $\mathfrak{R}_{\sigma}, \mathfrak{R}_{\tau}$ be compatible sets of positive roots for $G_{\sigma}$ and $G_{\tau}$. Then
\[
\g_{\sigma} /\g_{\tau} = \bigoplus_{\alpha \in \mathfrak{R}_{\sigma}\setminus \mathfrak{R}_{\tau}} \C_{\alpha}
\]
and 
\[
\mathrm{det}_{\C}(\nu_{\tau}^{\sigma}) = \bigotimes_{\alpha \in \mathfrak{R}_{\sigma}\setminus \mathfrak{R}_{\tau}} \C_{\alpha} = \C_{2(\rho_{\sigma} - \rho_{\tau})}. 
\]
\end{proof}
\end{lemma}

\begin{lemma}
There are canonical isomorphisms 
\[
\Psi_{\tau, \sigma}: S_{U_{\tau}}^{\mathrm{spin}}  \cong S_{U_{\sigma}}^{\mathrm{spin}} |_{U_{\tau}}, \hspace{5mm} \sigma \prec \tau
\]
and they automatically satisfy the cocycle condition. 
\begin{proof}
By the above lemma, we have that
\[
S_{Y_{\sigma}}|_{Y_{\tau}} \cong S_{(L\g)_{\sigma}/(L\g)_{\tau}}^{*} \otimes \C_{(\rho_{\sigma} - \rho_{\tau})}  \otimes S_{Y_{\tau}}.
\]
The claim follows form the construction. 
\end{proof}
\end{lemma}

\begin{definition}
We define the \emph{twisted spinor bundle} $S^{\mathrm{spin}}$ to be the $G$-equivariant bundle of  $\Z_{2}$-graded Hilbert space over $M$ with the property that
\[
S^{\mathrm{spin}}|_{U_{\sigma}} \cong S_{U_{\sigma}}^{\mathrm{spin}},
\] 
for all vertexes $\sigma$ of $\kA$. 
\end{definition}
The twisted spinor bundle is determined by the choice of almost complex structures on subsets $\{Y_{\sigma}\}$. By Lemma \ref{homotopy}, all the choices are homotopic. Hence, the twisted spinor bundle $S^{\mathrm{spin}}$ is unique up to homotopy.

\subsection{Construction of the twisted pre-quantum bundle}
Let us recall that the pre-quantum line bundle of a symplectic manifold is traditionally defined to be a line bundle whose first Chern class is an integral lift of the symplectic 2-form.
\begin{definition}\label{pre-quantum LG}
We say that a Hamiltonian $LG$-space $\M$ is \emph{pre-quantizable at level $k$} $(k >0)$ if there exists a $\widehat{LG}$-equivariant line bundle 
\[
\mathcal{L} \to \mathcal{M}
\]  
such that the central circle acts with weight $k$ and the first Chern class $c_{1}(\mathcal{L})$ equals to the symplectic 2-form on $\mathcal{M}$.
\end{definition}
Because the pre-quantum line bundle $\mathcal{L}$ is $\widehat{LG}$-equivariant instead of $LG$-equivariant, it might not descend to an actual line bundle on its corresponding q-Hamiltonian $G$-space $M$. 
\begin{remark}\label{pre-quantum q}
The 2-form $\omega$ for a q-Hamiltonian $G$-space $M$ is not closed in general. Instead the condition $d\omega = \phi^{*}\chi$ and the fact that$\chi$ is a closed 3-form imply that the pair $(\omega, \chi)$ defines a cocyle for the relative de Rham theory (see \cite[Appendix B]{MR2989614} for a reference). We denote by $[(\omega, \chi)] \in H^{3}(\phi, \R)$ its cohomology class. We say that a q-Hamiltonian $G$-space $(M, \omega, \phi)$ is \emph{pre-quantizable at level $k$} if $k \cdot [(\omega, \chi)]$ is integral. By the 1-1 correspondence between q-Hamiltonian $G$-spaces and Hamiltonian $LG$-spaces, their pre-quantum conditions are equivalent. 
\end{remark}

Let $Y_\sigma$ be the cross-section defined before. We identify it as a subset in $\M$.  If $\M$ has a pre-quantum line bundle $\cL$ at level $k$, then there exists a $\widehat{(LG)}_{\sigma}$-equivariant line bundle obtained by restriction
\[
L_{Y_{\sigma}} = \mathcal{L}|_{Y_{\sigma}} \to Y_\sigma,
\] 
on which the central circle acts with weight $k$. The collection of line bundles $\{L_{Y_{\sigma}}\}$ satisfy a compatibility condition in the sense that
\begin{equation}
\label{normal-2}
(LG)_{\sigma} \times_{(LG)_{\tau}}L_{Y_{\tau}} \cong L_{Y_{\sigma}}|_{Y_{\tau}^{\sigma}}, \hspace{5mm} \sigma \prec \tau. 
\end{equation}
where 
\[
Y_{\tau}^{\sigma} =(LG)_{\sigma} \times_{(LG)_{\tau}} Y_{\tau} 
\]
is a $(LG)_{\sigma}$-invariant open subset of $Y_{\sigma}$.

For any irreducible positive energy representation $V_{\lambda}$ at level $k$, we denote $V_{\lambda}^*$ its dual. Comparing to (\ref{spinor-2}), we define a bundle of Hilbert space on 
\[
U_\sigma = G\times_{G_\sigma} Y_\sigma
\] 
by
\begin{equation}
\label{pre-q-2}
S^{\mathrm{pre}}_{\lambda,U_{\sigma}} = G \times_{G_{\sigma}} (V_{\lambda}^* \otimes L_{Y_{\sigma}}).
\end{equation}
Here  the central circle of $\widehat{(LG)}_\sigma$  acts trivially on the tensor product $V_{\lambda}^* \otimes L_{Y_{\sigma}}$, and $G_\sigma$ acts factoring through $G_{\sigma} \cong (LG)_\sigma$.

By the compatibility condition (\ref{normal-2}), there are canonical isomorphisms 
\[
\Psi_{\tau, \sigma}: S^{\mathrm{pre}}_{\lambda,U_{\tau}}  \cong S^{\mathrm{pre}}_{\lambda,U_{\sigma}} |_{U_{\tau}}, \hspace{5mm} \sigma \prec \tau,
\]
satisfying the cocycle condition. Thus, they can be glued together. We define $S^{\mathrm{pre}}_{\lambda}$  the unique $G$-equivariant bundle of Hilbert space over $M$
with the property that 
\[
S^{\mathrm{pre}}_{\lambda}|_{U_{\sigma}} =S^{\mathrm{pre}}_{\lambda,U_{\sigma}}.
\] 

\begin{definition}
We define the \emph{twisted pre-quantum bundle} by 
\begin{equation}
\label{twisted pre-quantum}
S^{\mathrm{pre}}:= \bigoplus_{\lambda \in P_{k,+}} V_{\boldsymbol{\lambda}} \otimes S^{\mathrm{pre}}_{\lambda}. 
\end{equation}
It is a $\widehat{LG} \times G$-equivariant bundle of Hilbert space over the q-Hamiltonian $G$-space $M$. 
\end{definition}

\begin{remark}
There is a global construction of the twisted pre-quantum bundle.  Let us introduce  a $\widehat{LG} \times \widehat{LG}$-space:
\[
H_{\mathrm{wzw}, k} = \bigoplus_{\lambda \in P_{k,+}} V_{\lambda} \otimes V_{\lambda}^{*}.
\]
This is the so-called Hilbert space of the Wess-Zumino-Witten model (see \cite{Gawedzki00} for a reference). The aim of this paper is not to justify this choice of the Hilbert space, but  morally one can consider it as the analog of $L^2(G)$ for loop groups, in the spirit of the Peter-Weyl decomposition of $L^2(G)$. The tensor product
\[
(\M \times H_{\mathrm{wzw}, k}) \otimes \cL = \bigoplus_{\lambda \in P_{k,+}} V_{\lambda} \otimes (V_{\lambda}^{*} \otimes \cL) \to \M
\]
is a $\widehat{LG} \times LG$ equivariant bundle. We obtain a $\widehat{LG}\times G$-equivariant bundle of Hilbert space by taking its $\Omega G$-invariant part 
\[
S^{\mathrm{pre}} = V_{\lambda} \otimes [ V^{*}_{\boldsymbol{\lambda}} \otimes \mathcal{L}]^{\Omega G}  \to M = \M/ \Omega G.
\]
\end{remark}

\section{Dirac operators on q-Hamiltonian $G$-spaces}
With the twisted spinor bundle and twisted pre-quantum bundle defined in last section, we now proceed to construct Hilbert spaces and Dirac operators. We keep the same notations as in last section.  Let $G$ be a compact, simple and simply connected Lie group, and $M$ a pre-quantizable q-Hamiltonian $G$-space at level $k$.

\subsection{Dirac operators on cross-sections}
The idea of constructing Dirac operators is that we first define Dirac operators on local cross-sections, and then patch them together using partition of unity. 

Let $\{ Y_\tau \}$ be the collection of cross-sections defined in last section and open subsets $U_\tau = G \times_{G_\tau} Y_\tau$. 
Fixing an irreducible positive energy $LG$-representation $V_\lambda$, we define
\[
\big[\Gamma_{c}^\infty(U_{\tau}, S^{\mathrm{spin}}  \otimes S^{\mathrm{pre}}_{\lambda}) \big]^{G}
\]
the space of $G$-invariant, smooth sections of $S^{\mathrm{spin}} \otimes S^{\mathrm{pre}}_{\lambda}$ with compact support in $U_\tau$, with norm given by 
\[
\|s\|^2 := \int_{U_\tau} \langle s(m), s(m) \rangle dm.
\]
\begin{lemma}
\label{cross-1}
We have that
\begin{equation}
\begin{aligned}
&\big[\Gamma_{c}^\infty(U_{\tau}, S^{\mathrm{spin}}  \otimes S^{\mathrm{pre}}_{\lambda}) \big]^{G}\\
 \cong & \big[V_{\lambda}^{*} \otimes S^{*}_{L\g/(L\g)_{\tau}} \otimes \C_{(\rho_{G} -\rho_\tau, h^\vee)} \otimes  \Gamma_{c}^\infty(Y_\tau, S_{Y_\tau} \otimes L_{Y_\tau})\big]^{G_{\tau}}.
\end{aligned}
\end{equation}
\begin{proof}
The assertion follows immediately from the isomorphisms (\ref{spinor-2}) and (\ref{pre-q-2}). 
\end{proof} 
\end{lemma}

Let 
\[
D_{\mathrm{alg}} \in U(L\g) \otimes \mathrm{Cliff}(L\g/(L\g)_\tau)
\] 
be the cubic Dirac operator acting on $V_{\lambda}^{*} \otimes S_{L\g/(L\g)_\tau}^{*}$ and $D_{\mathrm{geo}}$ the equivariant  geometric Spin$^{c}$-Dirac operator on  $\Gamma_{c}^\infty(Y_\tau, S_{Y_\tau} \otimes L_{Y_\tau})$. Here we choose $D_{\mathrm{geo}}$ so that it is symmetric. Since the sum 
\[
D_{\mathrm{alg}} \otimes 1 + 1\otimes D_{\mathrm{geo}}
\]
is equivariant, it descends to the $G_\tau$-invariant part. That is, we obtain a collection of operators
\[
\big[D_{ \mathrm{alg}} \otimes 1 + 1\otimes D_{\mathrm{geo}}\big]^{G_\tau}
\] 
on
\[
 \big[V_{\lambda}^{*} \otimes S^{*}_{\p} \otimes \C_{(\rho_{G} -\rho_\tau, h^\vee)}  \otimes  \Gamma_{c}^\infty(Y_\tau, S_{Y_\tau} \otimes L_{Y_\tau})\big]^{G_\tau},
\]
and thus an operator on 
\[
\big[\Gamma_{c}^\infty(U_{\tau}, S^{\mathrm{spin}}  \otimes S^{\mathrm{pre}}_{\lambda}) \big]^{G}.
\]
We denote it by $D_{U_\tau}$. By definition, every $D_{U_\tau}$ is an unbounded operator. 
Since $D_{\mathrm{alg}}$ is defined only on vectors with finite energy,  the domain of $D_{U_\tau}$ is  given by 
\[
\big[(V_{\lambda}^{*} \otimes S^{*}_{L\g/(L\g)_{\tau}})^{\mathrm{fin}} \otimes \C_{(\rho_{G} -\rho_\tau, h^\vee)} \otimes  \Gamma_{c}^\infty(Y_\tau, S_{Y_\tau} \otimes L_{Y_\tau})\big]^{G_{\tau}},
\]
which is a dense subspace. 

Suppose that $\sigma_1, \sigma_2$ are two vertexes of $\kA$ and $\tau$ an open face of $\kA$ so that $ \sigma_1 \preceq \tau$ and $ \sigma_2 \preceq \tau$. In particular, one has that
\[
U_\tau \subseteq U_{\sigma_1}, \hspace{5mm} U_\tau \subseteq U_{\sigma_2}. 
\]

\begin{proposition}
\label{overlap}
If we restrict to 
\[
\big[\Gamma_{c}^\infty(U_{\tau}, S^{\mathrm{spin}}  \otimes S^{\mathrm{pre}}_{\lambda}) \big]^{G},
\] 
the difference between $D_{U_{\sigma_1}}$ and $D_{U_{\sigma_2}}$ is a bounded operator.
\begin{proof}
Since $U_\tau = G \times_{G_\tau} Y_\tau$, the Lemma \ref{cross-1} gives us an isomorphism:
\begin{equation} \label{iso-2}
\big[\Gamma_{c}^\infty(U_{\tau}, S^{\mathrm{spin}}  \otimes S^{\mathrm{pre}}_{\lambda}) \big]^{G} \cong  \big[  V_\lambda^*\otimes S^{*}_{L\g/(L\g)_{\tau}} \otimes \C_{(\rho_{G} -\rho_{\tau}, h^{\vee})} \otimes  \Gamma_{c}^\infty(Y_{\tau}, S_{Y_{\tau}} \otimes L_{Y_{\tau}})\big]^{G_{\tau}}.
\end{equation}
The Dirac operator $D_{U_\tau}$ is defined as the combination of a cubic Dirac operator on 
\[
V_\lambda^*\otimes S^{*}_{L\g/(L\g)_{\tau}} \cong  \cong V_\lambda^* \otimes S^{*}_{L\g/(L\g)_{\sigma_1}} \otimes S^*_{\g_{\sigma_1}/\g_\tau}
\]
and a Spin$^c$-Dirac operator on 
\[
  \Gamma_{c}^\infty(Y_{\tau}, S_{Y_{\tau}} \otimes L_{Y_{\tau}}).
\]

On the other hand, we have that
\[
U_{\tau} \cong G \times_{G_{\sigma_1}} Y_{\tau}^{\sigma_1}, \hspace{5mm} Y_{\tau}^{\sigma_1}= G_{\sigma_1} \times_{G_{\tau}} Y_{\tau}.
\]
Applying Lemma \ref{cross-1} again, 
\begin{equation}\label{iso-1}
\big[\Gamma_{c}^\infty(U_{\tau}, S^{\mathrm{spin}}  \otimes S^{\mathrm{pre}}_{\lambda}) \big]^{G}
\cong  \big[ V_\lambda^* \otimes S^{*}_{L\g/(L\g)_{\sigma_1}} \otimes \C_{(\rho -\rho_{\sigma_1}, h^{\vee})} \otimes  \Gamma_{c}^\infty(Y_{\tau}^{\sigma_1},S_{Y_{\sigma_1}} \otimes L_{Y_{\sigma_1}} )\big]^{G_{\sigma_1}}.
\end{equation}
Under the isomorphism (\ref{iso-1}), the operator $D_{U_{\sigma_1}}$ decomposes into the sum of a cubic Dirac operator on 
\[
V_\lambda^* \otimes S^{*}_{L\g/(L\g)_{\sigma_1}}
\]
and a Spin$^c$-Dirac operator on 
\[
 \Gamma_{c}^\infty(Y_{\tau}^{\sigma_1},S_{Y_{\sigma_1}} \otimes L_{Y_{\sigma_1}} ).
\] 
It follows immediately that the two operators  $D_{U_\tau}$ and $D_{U_{\sigma_1}}$ are identical on
\[
V_\lambda^* \otimes S^{*}_{L\g/(L\g)_{\sigma_1}}.
\] 

Note that the spinor bundle
\[
S_{Y_{\sigma_1}}|_{Y_{\tau}^{\sigma_1}} \cong  G_{\sigma_1} \times_{G_{\tau}} (S_{Y_{\tau}}  \otimes S^{*}_{\g_{\sigma_1}/ \g_{\tau}} \otimes \C_{\rho_{\sigma_1}-\rho_{\tau}}),
\]
and the pre-quantum line bundle
\[
L_{Y_{\sigma_1}}|_{Y_{\tau}^{\sigma_1}} \cong  G_{\sigma_1} \times_{G_{\tau}} L_\tau.
\]
Therefore,  the  space 
\[
\Gamma_{c}^\infty(Y_{\tau}^{\sigma_1},S_{Y_{\sigma_1}} \otimes L_{Y_{\sigma_1}} )
\] 
decomposes into:
\[
\big[ \big(C^\infty(G_{\sigma_1}) \otimes S^{*}_{\g_{\sigma_1}/ \g_{\tau}}\big) \otimes \C_{\rho_{\sigma_1}-\rho_{\tau}} \otimes \Gamma_{c}^\infty(Y_{\tau}, S_{Y_{\tau}} \otimes  L_{Y_{\tau}} )\big]^{G_{\tau}}. 
\]
By definition, both $D_{U_\tau}$ and $D_{U_{\sigma_1}}$ act as Spin$^{c}$-Dirac operators on 
\[
\Gamma_{c}^\infty(Y_{\tau}, S_{Y_{\tau}} \otimes  L_{Y_{\tau}} );
\] 
while on the factor
\[
C^\infty(G_{\sigma_1}) \otimes S^{*}_{\g_{\sigma_1}/ \g_{\tau}},
\] 
$D_{U_\tau}$ acts as the cubic Dirac operator, and $D_{U_{\sigma_1}}$ acts as the Spin$^{c}$-Dirac operator.  By Lemma \ref{geo alg}, their difference is bounded. 

We just show that the difference between $D_{U_{\sigma_1}}$ and $D_{U_{\tau}}$ is a bounded operator on the overlap $U_{\tau}$. Similarly one can show that the difference between $D_{U_{\sigma_2}}$ and $D_{U_{\tau}}$ is bounded as well. This completes the proof. 
\end{proof}
\end{proposition}

\subsection{Construction of Dirac operator and main theorem}
For a fixed irreducible positive energy representation $V_\lambda$ at level $k$, we define 
\[ 
\mathcal{H}_{\lambda} =  \big[L^2(M, S^{\mathrm{spin}}\otimes S^{\mathrm{pre}}_\lambda)\big]^{G}.
\]
The $\Z_2$-grading on the twisted spinor bundle $S^{\mathrm{spin}}$ equips $\cH_\lambda$ with a $\Z_2$-grading. That is 
\[
\cH_\lambda = \cH_\lambda^+ \oplus \cH_\lambda^-.
\]

Select a $G$-invariant, smooth partition of unity $\{f_{\sigma}^{2}\}$ which is subordinate to the cover 
\[
\{U_\sigma = G \times_{G_\sigma} Y_\sigma \}_{\sigma, \mathrm{dim}\sigma = 0}.
\] 
We define a Dirac operator $\cD$ on $\cH$ by the formula:
\[
\mathcal{D}= \sum_{\sigma, \mathrm{dim}\sigma = 0} f_{\sigma} \cdot  D_{U_\sigma} \cdot f_{\sigma}.
\]

\begin{proposition} \label{same up to bounded}
The Dirac operator $\cD$ doesn't depend on the choice of $\{f_\sigma^2\}$ up to homotopy.  
\begin{proof}
Let $\{f_{\sigma}'^{2}\}$ be another partition of unity. 
We define a Dirac operator $\cD'$ on $\cH_\lambda$ by
\[
\mathcal{D}'= \sum_{\sigma, \mathrm{dim}\sigma = 0} f'_{\sigma} \cdot  D_{U_\sigma} \cdot f'_{\sigma}.
\]
It suffices to show that 
\[
\mathcal{D}' - \mathcal{D} \in \mathbb{B}(\mathcal{H}_{\lambda}). 
\]

We compute that
\begin{equation}\label{step-1}
\begin{aligned}
\mathcal{D}&= \sum_{\sigma} f_{\sigma} \cdot  D_{U_\sigma} \cdot f_{\sigma} = \sum_{\sigma, \tau} f_{\sigma} \cdot  D_{U_\sigma}\cdot f'^2_\tau \cdot f_{\sigma}\\
& =  \sum_{\sigma, \tau} \big( f'_\tau \cdot f_{\sigma} \cdot  D_{U_\sigma} \cdot f_{\sigma} \cdot f'_\tau + f_{\sigma} \cdot  [D_{U_\sigma}, f'_\tau] \cdot f'_\tau \cdot f_{\sigma}\big)\\
& =  \sum_{\sigma, \tau} \big( f'_\tau  \cdot  D_{U_\sigma} \cdot f^2_{\sigma} \cdot f'_\tau  + f_{\sigma} \cdot  [D_{U_\sigma}, f'_\tau] \cdot f'_\tau \cdot f_{\sigma} + f'_{\tau} \cdot  [D_{U_\sigma}, f_\sigma] \cdot f'_\tau \cdot f_{\sigma}\big).\\
\end{aligned}
\end{equation}
Since the functions $f_\sigma$ are $G$-invariant, they commute with the cubic Dirac operator $D_{\mathrm{alg}}$. Their commutators with Spin$^c$-Dirac operators:
\[
c(df_\sigma) = [f_\sigma, D_{\mathrm{geo}}]
\]
are all bounded. It shows that 
\begin{equation}\label{step-2}
 [D_{U_\sigma}, f'_\tau],\hspace{5mm}   [D_{U_\sigma}, f_\sigma]  \in \mathbb{B}(\mathcal{H}_{\lambda}). 
\end{equation}

In addition, by Proposition \ref{overlap}, 
\[
 f'_\tau  \cdot ( D_{U_\sigma} - D_{U_\tau}) \cdot f^2_{\sigma} \cdot f'_\tau  
\]
are bounded for all $\sigma, \tau$. Therefore, 
\begin{equation}\label{step-3}
\begin{aligned}
&\sum_{\sigma, \tau} f'_\tau  \cdot  D_{U_\sigma} \cdot f^2_{\sigma} \cdot f'_\tau\\
&= \sum_{\tau, \sigma} f'_\tau  \cdot  D_{U_\tau}  \cdot f^2_\sigma \cdot f'_\tau + \sum_{\tau, \sigma}f'_\tau  \cdot ( D_{U_\sigma} - D_{U_\tau}) \cdot f^2_{\sigma} \cdot f'_\tau \\
& = \cD' + \sum_{\tau, \sigma}f'_\tau  \cdot ( D_{U_\sigma} - D_{U_\tau}) \cdot f^2_{\sigma} \cdot f'_\tau
\end{aligned}
\end{equation}
The proposition follows from (\ref{step-1}),  (\ref{step-2}) and  (\ref{step-3}). 
\end{proof}
\end{proposition}

\begin{proposition}\label{self-adjoint}
The Dirac operator $\cD$ is essentially self-adjoint. 
\begin{proof}
It is enough to show that each $f_\sigma \cdot D_{U_\sigma} \cdot f_\sigma$ is essentially self-adjoint. 
Let us write
\[
D_{U_\sigma} = \big[D_{ \mathrm{alg}} \otimes 1 + 1\otimes D_{\mathrm{geo}}\big]^{G_\sigma}
\] 
The algebraic part $D_{ \mathrm{alg}}$ is certainly self-adjoint, and the geometric part, $D_{\mathrm{geo}}$ is a symmetric Spin$^c$-Dirac operator on $Y_\sigma$. Recall that every symmetric Dirac operator on a complete manifold is essentially self-adjoint. But it doesn't apply directly to our case since $Y_\sigma$ is not complete. Nevertheless, we can get around this by the following trick.

Since the function $f_\sigma$ has compact support in $Y_\sigma$, we can find a smaller  subset $\widetilde{Y}_{\sigma}$ so that
\[
\mathrm{Supp} f_\sigma \subset \widetilde{Y}_{\sigma} \subset Y_\sigma
\]
and the closure of $\widetilde{Y}_{\sigma}$ is contained in $Y_\sigma$. We denote by $g$ the metric on $Y_\sigma$ and $\chi$ a positive function on $Y_\sigma$ so that $\chi|_{\widetilde{Y}_{\sigma}} \equiv 1$ and $\chi(m)$ tends to infinity as $m$ tends to the boundary of $Y_\sigma$.  Under the rescaled metric
\[
g^\chi(\cdot, \cdot) = \chi^2 \cdot g(\cdot, \cdot),
\]
the manifold $Y_\sigma$ becomes complete. 

Let $\widetilde D_{\mathrm{geo}}$ be an essentially self-adjoint Spin$^c$-Dirac operator  on the complete manifold $Y_\sigma$. Consider a new Dirac operator
\[
\widetilde D_{U_\sigma} := [D_{\mathrm{alg}} \otimes 1 + 1\otimes \widetilde D_{\mathrm{geo}}]^{G_{\sigma}},
\]
which is essentially self-adjoint. Since the metric remains the same within $\widetilde{Y}_\sigma$, we have that 
\[
\widetilde D_{U_\sigma} s = D_{U_\sigma} s
\] 
for all $s$ with $\mathrm{Supp} s \subseteq \widetilde{Y}_\sigma$. Thus, 
\[
f_\sigma \cdot \widetilde D_{U_\sigma} \cdot f_\sigma
\]
can be viewed as an operator on $\cH_\lambda$, and
\[
f_\sigma \cdot D_{U_\sigma} \cdot f_\sigma = f_\sigma \cdot \widetilde D_{U_\sigma} \cdot f_\sigma
\]
is essentially self-adjoint.

\end{proof}
\end{proposition}

We define an operator by functional calculus, 
\[
\cF = \frac{\cD}{\sqrt{1 + \cD^2}}.
\]
Because $\cD$ has a dense domain in $\cH_\lambda$, $\cF$ extends to a self-adjoint bounded operator on $\cH_\lambda$, which anti-commutes with the $\Z_2$-grading on $\cH_\lambda$.

\begin{theorem}\label{main theorem}
The bounded operator $\mathcal{F}$ is Fredholm on $\mathcal{H}_{\lambda}$. We define its Fredholm index by 
\[
\mathrm{Ind}(\cD)_{\lambda}:=  \mathrm{ker}(\mathcal{F}) \cap \mathcal{H}^{+}_{\lambda} - \mathrm{ker}(\mathcal{F}) \cap \mathcal{H}^{-}_{\lambda} \in \Z. 
\]
\end{theorem}
By Proposition \ref{same up to bounded},  the index is independent of the choice of partition of unity. Let us define a Hilbert space
\[
\mathcal{H} = \big[L^2(M, S^{\mathrm{spin}} \otimes S^{\mathrm{pre}})\big]^{G} \cong  \bigoplus_{\lambda \in P_{k,+}} V_{\lambda} \otimes \mathcal{H}_{\lambda}, 
\]
and a Dirac operator
\[
\cD_M = \bigoplus_{\lambda \in P_{k,+}} 1 \otimes \cD|_{\cH_\lambda}. 
\]
\begin{definition}
Let $(M, \omega)$ be a q-Hamiltonian $G$-space. If it is pre-quantizable at level $k$, we define its \emph{quantization} 
\[
Q(M) =\mathrm{Ind}(\cD_M) =  \bigoplus_{\lambda \in P_{k,+}} \mathrm{Ind}(\cD)_{\lambda}  \cdot  V_{\lambda} \in R_{k}(LG).
\]
\end{definition}
This generalizes Bott's Spin$^c$-quantization for Hamiltonian $G$-spaces to q-Hamiltonian $G$-spaces.

\begin{remark}
In \cite{MR2989614} Meinrenken develops a quantization from pre-quantized q-Hamiltonian $G$-spaces to the equivariant twisted $K$-homology of $G$ using push-forward maps.  By the work of Freed-Hopkins-Teleman \cite{Freed13} the equivariant twisted $K$-homology of $G$ at level $k$ is isomorphic to fusion ring of loop group $R_{k}(LG)$.  
\end{remark}

\begin{example}
Let 
\[
\mathcal{C} = G \cdot \exp(\xi) \cong G/H, \hspace{5mm} \xi \in \kA
\] 
be a conjugacy class. We assume that it is pre-quantizable at level $k$. By Remark \ref{pre-quantum q},   
\[
(\eta, k) = (k \cdot \xi, k) \in \widehat{\kt}^*_+
\] 
is an integral weight. By definition, the twisted spinor bundle and twisted pre-quantum bundle are given by 
\[
S^{\mathrm{spin}} = G \times_{H} (S_{L\g, \h}^* \otimes \C_{(\rho_{G} - \rho_{H}, h^{\vee})}), \hspace{5mm} S^{\mathrm{pre}} = G\times_{H}(H_{\mathrm{wzw}, k} \otimes \C_{(\eta, k)}).
\]
The Hilbert space 
\[
\big[L^2(\mathcal{C}, S^{\mathrm{spin}}\otimes S^{\mathrm{pre}})\big]^G  \cong
\bigoplus_{\lambda \in P_{k,+}} V_{\boldsymbol{\lambda}} \otimes \big(V_{\boldsymbol{\lambda}}^{*} \otimes S_{L\g/\h}^{*} \otimes \C_{(\rho_{G}-\rho_{H} +\eta, k +h^{\vee})} \big)^{H}.
\]
By Theorem \ref{square operator-4}, one can calculate that
\[
Q(\mathcal{C})= \mathrm{ind}(\mathcal{D}_M) = V_{\eta} \in R_{k}(LG).
\]
This is an algebraic version  of the Borel-Weil construction for loop groups \cite{PS86}. 
\end{example}

\subsection{Proof of Theorem \ref{main theorem}}
We will prove the main theorem in this subsection. Let us begin with a lemma. 
\begin{lemma}
\label{bounded perturb}
Let $\cD,\cD'$ be two self-adjoint unbounded operators on a Hilbert space $\cH$ such that $\cD - \cD' \in  \mathbb{B}(\cH)$. For $\alpha = 0,1$, if 
\[
\cD^\alpha \cdot (1+ \cD^2)^{-1} \in \mathbb{K}(\cH), 
\]
then 
\[
\cD'^\alpha \cdot (1+ \cD'^2)^{-1} \in \mathbb{K}(\cH).
\]
\begin{proof}
Suppose that $\cD' = \cD + B$ with  $B \in  \mathbb{B}(\cH)$. By straightforward calculation,
\begin{equation}\label{op-2}
\begin{aligned}
& (1+ \cD^2)^{-1} -  (1+ \cD'^2)^{-1} \\
& =(1+ \cD^2)^{-1} \cdot(\cD'^2 - \cD^2) \cdot (1+ \cD'^2)^{-1}\\
& = (1+ \cD^2)^{-1} \cdot(\cD \cdot B + B \cdot \cD') \cdot (1+ \cD'^2)^{-1}
\end{aligned}
\end{equation}
By the fact that the product of bounded operator and compact operator is again a compact operator and the assumption, we deduce that
\[
(1+ \cD'^2)^{-1}\in \mathbb{K}(\cH).
\]

We next consider
\begin{equation}\label{op-2}
\begin{aligned}
& \cD \cdot (1+ \cD^2)^{-1} -  \cD \cdot  (1+ \cD'^2)^{-1} \\
& = \cD \cdot (1+ \cD^2)^{-1} \cdot(\cD'^2 - \cD^2) \cdot (1+ \cD'^2)^{-1}\\
& = \cD \cdot (1+ \cD^2)^{-1} \cdot(\cD \cdot B + B \cdot \cD') \cdot (1+ \cD'^2)^{-1}
\end{aligned}
\end{equation}
For the same reason as above, it follows that
\[
\cD \cdot (1+ \cD'^2)^{-1}\in \mathbb{K}(\cH).
\]
Because $\cD' -\cD$ is bounded, we conclude that 
\[
\cD' \cdot (1+ \cD'^2)^{-1}\in \mathbb{K}(\cH).
\] 
\end{proof}
\end{lemma}

Fix a vertex $\sigma$ of $\kA$. Let us write
\[
D_{U_\sigma} = \big[D_{ \mathrm{alg}} \otimes 1 + 1\otimes D_{\mathrm{geo}}\big]^{G_\sigma},
\] 
acting on
\[
 \big[V_{\lambda}^{*} \otimes S^{*}_{L\g/(L\g)_{\sigma}} \otimes \C_{(\rho_{G} -\rho_\sigma, h^\vee)} \otimes  \Gamma_{c}^\infty(Y_\sigma, S_{Y_\sigma} \otimes L_{Y_\sigma})\big]^{G_{\sigma}}.
\]
By the trick used in Proposition \ref{self-adjoint}, we can assume that $Y_\sigma$ is complete and $D_{\mathrm{geo}}$ is an essentially self-adjoint Spin$^c$-Dirac operator on $Y_\sigma$. 

 \begin{lemma}\label{compact op}
For $\alpha = 0,1$, one has that
\[
f_\sigma \cdot  D^{\alpha}_{U_\sigma} \cdot (1+ D^{2}_{U_\sigma})^{-1}  \cdot f_\sigma  \in \mathbb{K}(\cH_\lambda). 
\] 
\begin{proof}
We decompose 
\[
L^2(Y_\sigma, S_{Y_\sigma} \otimes L_{Y_\sigma})
\] 
with respect to the $S^{1}_{\mathrm{rot}} \times \widehat{G}_\sigma$-action and denote by $M(\boldsymbol{\nu})$ the isotypic  component labeled by 
\[
\boldsymbol{\nu} = (n, \nu, k) \in \Z \times \Lambda^{*} \times \Z.
\] 
Since $f_{\sigma}$ has compact support, it follows from the Rellich's lemma that 
\[ f_\sigma \cdot D_{\mathrm{geo}}^{\alpha} \cdot  (1 + D_{\mathrm{geo}}^{2} )^{-1} \cdot f_{\sigma} 
\]
is a compact operator on 
\[
L^2(Y_\sigma, S_{Y_\sigma} \otimes L_{Y_\sigma}).
\] 
It implies that the norm of the restriction of 
\[
f_\sigma \cdot  D^{\alpha}_{U_\sigma} \cdot (1+ D^{2}_{U_\sigma})^{-1}  \cdot f_\sigma
\]
to 
\[
\big[ V_{\lambda}^*\otimes S^{*}_{L\g/(L\g)_{\sigma}} \otimes \C_{(\rho_{G} -\rho_{\sigma}, h^{\vee})} \otimes  M(\boldsymbol{\nu})\big]^{G_{\sigma}}.
\]
tends to 0 as $\|\boldsymbol{\nu}\|$ tends to infinity. Thus, it suffices to show that the operator
\[
f_\sigma \cdot  D^{\alpha}_{U_\sigma} \cdot (1+ D^{2}_{U_\sigma})^{-1}  \cdot f_\sigma
\]
is compact after restricted to a fixed component 
\[
\big[ V_{\lambda}^*\otimes S^{*}_{L\g/(L\g)_{\sigma}} \otimes \C_{(\rho_{G} -\rho_{\sigma}, h^{\vee})} \otimes  M(\boldsymbol{\nu})\big]^{G_{\sigma}}.
\]

By Theorem \ref{square operator-4},  the operator $D^2_{\mathrm{alg}} = D^2_{L\g,(L\g)_\sigma}$ acts on the above space by 
\[
 \|\boldsymbol{\lambda} + \boldsymbol{\rho}_{G}\|^{2} - \|\boldsymbol{\nu} + \boldsymbol{\rho}_{G}\|^{2}.
\]
By the formula in (\ref{weight product}), 
\[
 \|\boldsymbol{\lambda} + \boldsymbol{\rho}_{G}\|^{2} - \|\boldsymbol{\nu} + \boldsymbol{\rho}_{G}\|^{2} =  \|\boldsymbol{\lambda} + \boldsymbol{\rho}_{G}\|^{2} + 2 \cdot n \cdot (k+ h^\vee) - \|\nu + \rho_G\|^2.
\]
Since $ \|\boldsymbol{\lambda} + \boldsymbol{\rho}_{G}\|^{2}$ and $ \|\nu + \rho_g\|^2$ are fixed constants, the operator $D^2_{\mathrm{alg}}$ tends to infinity as the energy level $n$ goes to infinity. This proves the assertion. 
\end{proof}
\end{lemma}

\begin{lemma}\label{key lemma}
One has that
\[
f_\sigma \cdot (1+\mathcal{D}^{2})^{-1}  \cdot f_\sigma, \hspace{5mm} [f_\sigma,  (1+\mathcal{D}^{2})^{-1}] \cdot f_\sigma \in \mathbb{K}(\mathcal{H}_{\lambda})
\]
\begin{proof}
As shown in Proposition \ref{same up to bounded}, there exists a bounded operator $B$ so that
\[
\cD \cdot s = D_{U_\sigma}\cdot s + B \cdot s.
\]
for any $s \in \cH_\lambda$ with $\mathrm{Supp} s \subset U_\sigma$. By the choice of function $f_\sigma$,  the support of $f_\sigma \cdot s$ is automatically contained in $U_\sigma$ for all $s \in \cH_\lambda$. One can verify that
\[
f_\sigma \cdot \cD^\alpha \cdot  (1+\mathcal{D}^{2})^{-1}  \cdot f_\sigma = f_\sigma \cdot (D_{U_\sigma} + B)^\alpha \cdot \big(1+(D_{U_\sigma} + B)^{2} \big)^{-1}  \cdot f_\sigma.
\]
By Lemma \ref{bounded perturb} and \ref{compact op}, we conclude that for $\alpha=0,1$, 
\[
f_\sigma \cdot \cD^\alpha \cdot  (1+\mathcal{D}^{2})^{-1}  \cdot f_\sigma \in \mathbb{K}(\cH_\lambda).
\]

For the second half, we calculate that
\begin{equation} 
\begin{aligned}
&\big[f_\sigma,  (1+\mathcal{D}^{2})^{-1}\big]\cdot f_\sigma \\
= &(1+\mathcal{D}^{2})^{-1} \cdot [f_\sigma, \cD^2] \cdot(1+\mathcal{D}^{2})^{-1}\cdot f_\sigma\\
=& (1+\mathcal{D}^{2})^{-1} \cdot c(df_\sigma)\cdot \cD   \cdot(1+\mathcal{D}^{2})^{-1}\cdot f_\sigma\\
+& (1+\mathcal{D}^{2})^{-1} \cdot \cD \cdot c(df_\sigma) \cdot(1+\mathcal{D}^{2})^{-1}\cdot f_\sigma.
\end{aligned}
\end{equation}
We point out here that $c(df_\sigma)$ is a bounded operator with support  in $U_\sigma$. As in Lemma \ref{compact op}, one can similarly prove that
\[
 c(df_\sigma)\cdot \cD^\alpha   \cdot(1+\mathcal{D}^{2})^{-1}\cdot f_\sigma \in \mathbb{K}(\cH_\lambda). 
\]
This completes the proof. 
\end{proof}
\end{lemma}

Because the operator
\[
\cF = \frac{\cD}{\sqrt{1 + \cD^2}} \in \mathbb{B}(\cH_\lambda),
\]
we have that 
\begin{equation}
\begin{aligned}
 &1 - \cF^2 = (1+\mathcal{D}^{2})^{-1}\\
 = &\sum_{\mathrm{dim}\sigma = 0}  (1+\mathcal{D}^{2})^{-1} \cdot f^2_\sigma\\
 =&  \sum_{\mathrm{dim}\sigma = 0} \big( f_\sigma \cdot (1+\mathcal{D}^{2})^{-1} \cdot f_\sigma +  [ f_\sigma,  (1+\mathcal{D}^{2})^{-1}] \cdot f_\sigma \big). 
\end{aligned}
\end{equation}
By Lemma \ref{key lemma}, 
\[
 1 - \cF^2 \in \mathbb{K}(\mathcal{H}_{\lambda}),
\]
which implies that $\cF$ is a Fredholm operator on $\mathcal{H}_{\lambda}$.

\bibliographystyle{alpha}
\bibliography{mybib}

\end{document}